\documentclass[10pt,a4paper]{article}

\usepackage{cogsci}
\usepackage{pslatex}
\usepackage{apacite}
\usepackage[dvipdfmx]{graphicx}
\usepackage{color}

\usepackage{amsmath, amssymb}
\usepackage[all]{xy}

\newtheorem{thm}{Theorem}
\newtheorem{defi}[thm]{Definition}

\newcommand{\MF}[1]{\textcolor{black}{#1}}

\DeclareMathOperator{\dom}{dom}
\DeclareMathOperator{\cod}{cod}

\title{Meanings, Metaphors, and Morphisms: \\Theory of Indeterminate Natural Transformation (TINT)}
 
\author{{\large \bf Miho Fuyama (miho02@sj9.so-net.ne.jp)} \\
  Tamagawa University Brain Science Institute\\6-1-1 Tamagawagakuen Machida-shi, 
  Tokyo, 194-8610, Japan
  \AND {\large \bf Hayato Saigo (harmoniahayato@gmail.com)} \\
  Nagahama Institute of Bio-Science and Technology\\
   1266 Tamura-cho, Nagahama-shi, Shiga, 526-0829, Japan}

\begin{document}

\maketitle

\begin{abstract}
In the present paper\MF{,} we propose a new theory named ``\textit{Theory of indeterminate natural transformation (TINT)}'' 
to investigate the dynamical creation of meanings as association relationships between images, focusing on the metaphor comprehension as an example. 
TINT models the meaning creation as a kind of stochastic processes based on 
the mathematical structure defined by association relationships as morphisms in category theory,  
so as to represent the indeterminate nature of structure-structure interactions between the systems of the meanings of images. Such interactions are 
formulated in terms of so-called coslice categories and functors as structure\MF{-}preserving correspondence between them. 
The relationship between such functors \MF{is} ``indeterminate natural transformation'', the central notion in TINT, which models the creation of meanings 
in a precise manner. For instance, the process of metaphor comprehension is modeled by the construction of indeterminate natural transformation from a 
canonically defined functor which we call the base-of-metaphor functor. 
\end{abstract}

\section{1. Introduction to TINT}
Our minds are filled with various images, visual or auditory, concrete or abstract, verbal or non-verbal:
The image of a cup, snow, mothers, a girl called Naomi, red, cold, a series of sounds like``kiki", the God, Pegasus, the number 3, and so on. 
For each image, plenty of meanings are being created in our minds.  
This ``creation of meanings'' is one of the essential abilities of humankind, 
and also one of the fascinating topics for researchers interested in cognitive processes. 

In the present paper\MF{,} we investigate the process of creating meanings, especially the creations of meanings of novel metaphors
\footnote{In the present paper\MF{,} we do not distinguish similes from metaphor, for our focus is the general structures not dependent on such classification.} 
as typical examples, 
from the viewpoint that the meanings of an image are nothing but the relationships to other images. 
As a suitable framework \MF{for} the investigation\MF{,} we propose 
a new theory of metaphor, \textit{Theory of Indeterminate Natural Transformations} (\textit{TINT}). 

The aim of TINT is to describe and explain how meanings of the source and the target of a metaphor
\footnote{``Target'' and ``source'' are terminology of metaphor study. For example, in the case of ``A is B'' or ``A is like B'', 
A is a target and B is a source of this metaphor.}
 interact with each other and create novel meanings based on the structure of relationships among images. 
TINT is formulated on 
the fundamental concepts in category theory which \MF{is a branch of} mathematics to investigate various mathematical structures 
in terms of relationships called ``morphisms''(or ``\MF{arrow}''). More precisely, 
we extend the notion of category, functor and natural transformation to the context of a kind of stochastic process to model the ``indeterminate'' 
nature of the creation of meanings.

The contents of the present paper are as follows: In section 2, we briefly explain the intuitive explanation 
of TINT without using \MF{the} terminology of category theory.
Then we formulate the fundamental concept of TINT such as indeterminate natural transformation, after introducing and extending 
basic concepts in category theory in section 3. 
The central notion is ``indeterminate natural transformation''.
The following section 4 is devoted to give an example of the applications of TINT to metaphor comprehension. In the last section 
we discuss some essential problems of metaphor comprehension based on TINT and some prospects for future investigations.

\section{2. Explanation of TINT}
First of all,
we try to define a ``meaning'' of image.
In TINT,
a meaning of a certain image is defined as a whole associative relationships  from its image to other images.
Suppose some sort of event raise image $A$, and the image $A$ evokes another image $B$, $C$, and so on.
Then we define the meaning of $A$ as the whole association from image $A$ to image $B$, $C$ and some other images.
For example, when a image of ``Love'' evokes some images such as ``Warm'', ``Children'', ``Lover'', ``Suffer'', ``Passion'', ``Heart'', etc, we regard the meaning of ``Love'' as the whole association from ``Love'' to these ``Warm'', ``Children'', ``Lover'', ``Suffer'', ``Passion'', ``Heart'', etc.

According to this line of thinking, creating new meanings can be recognized as the creating new relationships between images.
Therefore, the simplest case of creating meaning is that one associative relationship from an image $A$ to an image $B$ are raising and this phenomena changes the meaning of $A$ and create new meanings of other images.

An example of this process is comprehension of novel metaphor,
 in which a new connection between a target and a source of the metaphor is excited.
For example, 
if the metaphor of ``Love is drinking water'' is new for you,
an associative relationship from an image of ``Love'' to an image of ``drinking water'' is made.
Then the meaning of ``Love'' is changed by its new relationship.
In this manner, we can recognize the process of comprehending novel metaphors as one of the most basic processes of creating new meanings.

Furthermore,
each image of the target and the source have had associations to other images before understanding  new metaphor.
Therefore, when comprehending the metaphor, 
we will not only make the associative relationship from ``Love'' to ``Drinking water'' but also correspond some part of associative relationships from ``Love'' to some part of associative relationships from ``Drinking water''.
These associating each structure fertilize the meaning of ``Love''; Love is like drinking water, so Love is necessary for live, feel good, harmful without cleanliness, transparent, cold, and so on.

We think it is the most important process of metaphor comprehension to
this constructing relationship between a part of associative structures of the target and the source.
TINT represents this structure-to-structure interaction as construction of \textit{indeterminate natural transformation}. The next section is 
devoted to formulate this central concept based on the fundamental notions in category thoery and its stochastic extension.

\section{3. Formulation of TINT}
\subsection{Category-theoretic concepts for TINT}

The fundamental concept of the categorical theory as referred to here are "category", "functor", "natural transformation"
 (for more detailed information, see \cite {Mac1998}.).

A category is, roughly speaking, a network formed by composable "morphisms" which intertwines "objects". 
It can be considered that the objects represent some "phenomena", and the morphisms \MF{represent
transformations or processes}  between those phenomena.

0. A \textbf {category} is a system consisting of \textbf{objects} and \textbf {morphisms}, satisfying the following four conditions.\\

I. Each morphism $f$ is associated with two objects $\dom(f)$ and $\cod(f)$, 
which are called \textbf{domain}, \textbf{codmain}. 

When $\dom(f)=X$ and $\dom(f)=Y$, we denote
\begin{align*}
f: X \longrightarrow Y
\end{align*}
or
\begin{align*}
X \xrightarrow[]{~f~} Y
\end{align*}
\textit{It is not necessary to limit the direction of the morphism from left to right, 
if it is convenient it is free to write from bottom to top, or from right to left, etc..}
A subsystem of the category build up with these morphisms and objects are called \textbf{diagrams}.\\

II. If there are two morphisms $f,g$ such that $\cod(f)=\dom(g)$, in other words,
\begin{align*}
Z \xleftarrow[]{~g~} Y \xleftarrow[]{~f~} X
\end{align*}
there is a unique morphism
\begin{align*}
Z \xleftarrow[]{~g \circ f~} X
\end{align*}
called the \textbf{composition} of $f,g$.\\

III. We assume so-called \textbf{associative law}. For the diagram
\begin{align*}
\xymatrix@C=.5cm@R=\halfrootthree cm{
W && Y \ar[dl]^g \ar[ll]_{h \circ g} \\
 & Z \ar[ul]^h && X \ar[ul]_f \ar[ll]^{g \circ f}
}
\end{align*}
we assume
\begin{align*}
(h \circ g) \circ f
= h \circ (g \circ f)
\end{align*}
As above, when all compositions of morphisms which \MF{have} the same codomain and domain are equal,  we call the diagram \textbf{commutative}.\\

IV. The last condition is \textbf{unit law}: For any object $X$ there exist a morphism 
$1_X : X \longrightarrow X$ called the \textbf{identity} of 
$X$ such that the diagram

\begin{align*}
\xymatrix@C=.5cm@R=\halfrootthree cm{
& X \ar[ld]_f \\
Y && X \ar[ll]_f \ar[ld]^f \ar[lu]_{1_X} \\
& Y \ar[lu]^{1_Y}
}
\end{align*}
is commutative for any $f: X \longrightarrow Y$. In other words, 
\begin{align*}
f \circ 1_X = f = 1_Y \circ f.
\end{align*}
By the natural correspondence from objects to their identities, we may ``identify'' the objects as identities. In other words, we may consider the objects are 
just the special morphisms. In the following sometimes we may adopt this viewpoint without notice.

\begin{defi}
A category is a system composed of two kinds of entities called objects and morphisms, which are interrelating through the notion of domain/codomain, 
equipped with composition and identity, satisfying associative law and unit law.

\end{defi}

\paragraph{Example: Category of images $\mathcal{C}$}
The objects of $\mathcal{C}$ are images (visual or auditory, concrete or abstract, verbal or nonverbal) and morphisms of $\mathcal{C}$ are associations between them.
\footnote{For simplicity, in the present paper\MF{,} we think the number of morphisms between two given object is at \MF{most one}. 
That is, the morphism is considered as ``the possibility 
of association'', though it is fruitful to relax this condition to the general case (i.e. the case with \MF{a} distinction between many kinds of association). }\\

\paragraph{Example: Category of meanings as coslice category} Let $A$ be an object of $\mathcal{C}$. 
The coslice category $A\backslash \mathcal{C}$ is defined as follows:
Objects of $A\backslash \mathcal{C}$ is morphism from $A$. The morphisms between two objects $f_1:A\longrightarrow X_1$ and $f_2:A\longrightarrow X_2$ 
is the triple of $(f_1,f_2,g)$ where 
$g:X_1\longrightarrow X_2$ and $f_2=g\circ f_1$. That is, the objects are morphisms from $A$ (Something ``for $A$'') and 
the morphisms are commutative triangle diagrams (relationships for $A$). 
When $\mathcal{C}$ is the category of images, we identify this category as the categories of meanings of $A$.\\

A functor is defined as a structure-preserving correspondence of two categories:

\begin{defi}

A correspondence $F$ from $\mathcal{C}$ to $\mathcal{D}$ which maps each object/morphism in $\mathcal{C}$ to corresponding object/morphism in $\mathcal{D}$ is called 
a \textbf{functor} if is satisfies the following conditions:

\begin{enumerate}
\item It maps $f:X \longrightarrow Y in $ $\mathcal{C}$ to $F(f):F(X) \longrightarrow F(Y)$ in $\mathcal{D}$.
\item $F(f \circ g) =F(f) \circ F(g)$ for any (composable) pair of $f,g$ in $\mathcal{C}$. 
\item For each $X$ in $\mathcal{C}$, $F(1_X) = 1_{F(X)}$.
\end{enumerate}
\end{defi}

In short, a functor is a correspondence which preserves diagrams, or equivalently, preseves categorical structure.

\MF{The functor} is a \MF{universal} concept. All the processes 
expressed by words such as recognition, representation, construction, modeling, theorization, etc. can be said to be the creation of functors.

\paragraph{Example (Metaphor as Functor)}
A conceptual metaphor itself can be considered as functors. For example, 
the metaphor from atomic worlds to celestial worlds is considered as some functor from the categories whose objects are protons, electrons\MF{,} 
and electromagnetic force etc. and morphisms are associations between them 
to those of stars, planets and gravitational forces etc. and associations between them.\\

The examples above are the ``completed metaphors''. The main topic in the present paper is, rather the \textit{construction} of such metaphors.
To model the construction itself it is important to focus on the ``Base-of-metaphor'' functor between coslice categories:

\paragraph{Example (Base-of-Metaphor functor between the categories of meanings)}
Let us consider a given morphism $f:A\longrightarrow B$ the association induced by 
a (new) metaphor expresseion as ``A is like B''. Then a functor $f\backslash \mathcal{C}:=(\cdot)\circ f: B\backslash \mathcal{C}\longrightarrow 
A\backslash \mathcal{C}$ (i.e. ``$(\cdot)$ for $B$'' for $A$ )is canonically defined. We call it the base-of-metaphor functor.

Our basic hypothesis is that the construction of \MF{the} meaning of a new functor (metaphor functor) from base-of-metaphor functor. 
Of course\MF{,} the base-of-metaphor functor itself is a functor and provides some indirect meanings like ``$X$ for $B$ for $A$'', our minds tend to construct 
other new association from it to some ``$Y$ for $A$'': ``$X$ is like $Y$ (for $A$)''. 
The construction should be considered as some ``morphism'' from the base-of-metaphor functor to some other new metaphor functor. Then, 
what is the morphism between functors? The answer is natural transformation:

\begin{defi}
Let $F,G$ be functors from category $\mathcal{C}$ to category $\mathcal{D}$, a correspondence $t$ is called a \textbf{natural transformation} 
from $F$ to $G$ if it satisfies the following conditions:
\begin{enumerate}
\item $t$ maps each object $X$ in $\mathcal{C}$ to corresponding morphism $t_X:F(X)  \longrightarrow G(X)$ in $\mathcal{D}$.
\item For any $f:X \longrightarrow Y$ in $\mathcal{C}$, 
\begin{align*}
t_Y \circ F(f) = G(f) \circ t_X.
\end{align*}
\end{enumerate}
\end{defi}

For the natural transformation we use the notation such as $t:F\Longrightarrow G$. The second condition above is depicted as below:

\begin{align*}
\xymatrix{
\ar@<-10pt>@{}[d]_~="t1" \ar@<10pt>@{}[r]^~="t3"
& Y & X \ar[l]_f \ar@<10pt>@{}[d]^~="t2" \\
F \ar@{=>}[d]_t & F(Y) \ar[d]_{t_Y} & F(X) \ar[d]^{t_X} \ar[l]_{F(f)} \\
G \ar@<-10pt>@{}[r]_~="t4" & G(Y) & G(X) \ar[l]^{G(f)}
\ar@{.} "t1";"t2"
\ar@{.} "t3";"t4"
}
\end{align*}
Upper-right part denotes the morphism in $\mathcal{C}$ and downer-left part the morphism in $\mathcal{D}$. 
The second condition in the definition of natural transformation 
means that the diagram above is commutative.\footnote{Note that in the setting of the present paper (i.e. there is at most one morphism between given two objects) 
this condition is automatically satisfied when the functors are 
well-defined, though in the future work we will focus on the general case.}

\MF{
In sum,
figure 1 shows the schematic diagram of this process from evoking morphism $f:A\longrightarrow B$ to construction of a new metaphor functor.}

\begin{figure}[h]
\begin{center}
\includegraphics[width = 6cm,height=4cm]{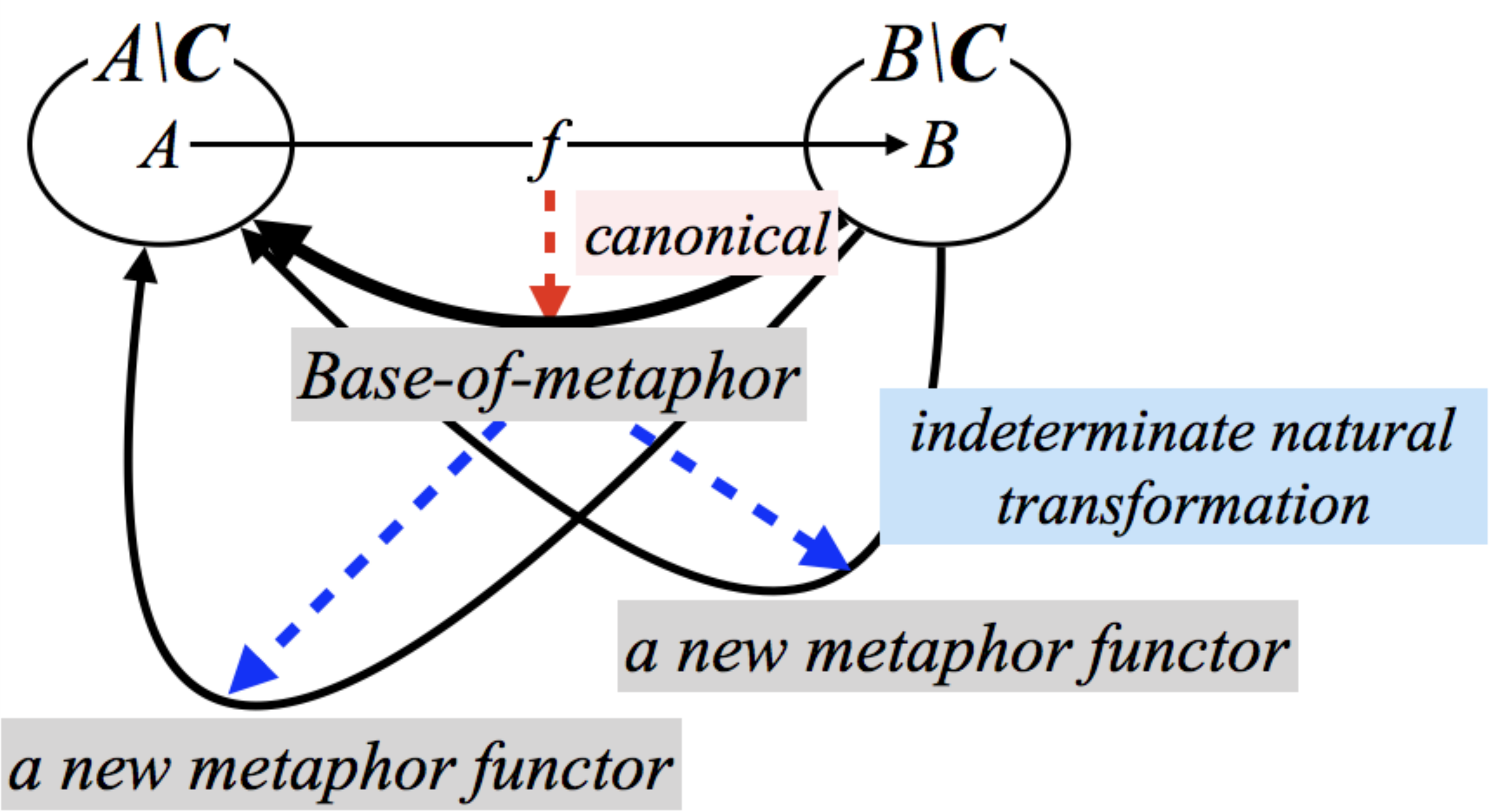}
\end{center}
\caption{\MF{Base-of-metaphor is canonically made from morphism $f$, and the metaphor functor is constructed based on this base-of-metaphor by searching natural transformation.}} 
\label{Fig_base} 
\end{figure}

We have introduced the basic concepts of category theory to formulate meanings and metaphors in mathematical terms. 
However, the most important topic in the present paper is 
not the completed metaphor functors and natural transformation, but the process to create them. 
As all of us understand in the daily life that the possible interpretation of a 
metaphor expression is by no means unique or completed, and that the process of creating \MF{a} new interpretation of a metaphor cannot be considered as \MF{the}
deterministic process as in classical mechanics. 
In the next subsection, we introduce the concept of indeterminate category, indeterminate functor\MF{,} and indeterminate natural \MF{transformation}
to treat with this kind of indeterminacy on the basis of category theory.
 
\subsection{Toward Indeterminacy}

The intuitive ideas of \MF{an} indeterminate category, indeterminate functor and indeterminate natural transformation themselves are quite simple. They are just 
fluctuating, growing and co-existing family of (subsystem of) categories, functors and natural transformation, constructed through indeterministic processes. To model this simple concept, 
we need to introduce some stochastic concepts over categories such as \MF{the} weight of morphisms, \MF{the} excitement of morphisms and excitation/relaxation processes. For simplicity, 
we treat the case that the morphism has two-level states: excited or relaxed. 

A stochastic category
 $(\mathcal{C},\mu,R)$ is the triple of 
a category $\mathcal{C}$, 
the probabilistic weight $\mu$ (without normalization) on morphisms, and the excitation/relaxation rules $R$ referring to $\mu$. 
We assume $R$ contains the following rules,

\begin{itemize}
 \item (Basic rule). 
 \begin{itemize}
 \item The morphism composed by two excited morphisms is excited.
 \item Identities of the domain and codomain of excited morphisms are excited. 
 \end{itemize}
\end{itemize}
which themselves are not depending on $\mu$. Based on the notion of stochastic category, a "stochastic process" of category $\mathcal{C}_{exc,t}$ 
as the (subsystem of) category of $C$ containing all the excited morphisms in the step $t$.
In the present paper we consider this $\mathcal{C}_{exc,t}$ as the mathematical realization of the concept of 
``indeterminate categories''.
\MF{An} indeterminate functor or indeterminate natural transformation id defined by (family of the subsystem of) functor and natural 
transformations between indeterminate categories, constructed in terms of stochastic categories.

It seems very interesting to study the concept and examples mathematically and application to various scientific branches. In the present paper, however, 
we focus on the modeling of the construction of meanings of a metaphor, i.e. TINT.
 
\subsection{Axioms of TINT}

Based on the all the arguments above, we propose (a sketch of) axioms of TINT as working hypothesis:

\begin{itemize}
 \item The system of all the images (visual or auditory, concrete or abstract, verbal or non-verbal) 
 and all the association between them can be modeled by a category $\mathcal{C}$. Stochastic aspect is provided by 
 a weight on $\mu$ and excitation/relaxation rule introduced below:
 \item The system of meanings of an image $A$ and relation between meanings can be modeled by coslice category and weight on it naturally induced by $\mu$:
 \item A Metaphor ``$A$ is like $B$'' excites the morphism $f:A\to B$\MF{.} \footnote{When the weight of the morphism, is (almost) zero, maybe ``creates'' 
 is more suitable than ``excites''.} It \MF{causes} the construction of natural transformations based on 
 the excitation/relaxation processes described below, which are nothing but the basic dynamics of development of coslice categories:

 \item Excitation processes (relatively fast) occurs under the stochastic rules as below:
 \begin{itemize}
 
 \item 0. (Basic rule). The morphism composed of two excited morphisms is excited. Identities of the domain and codomain of excited morphisms are excited. 
 \item 1. (Neighboring rule). The morphisms whose domain is the codomain of an excited morphism tend to be excited: The excitation probability is provided by $\mu$. 
 \item 2. (Fork rule). For the pairs of morphisms sharing the domains, morphisms between their codomains (which may be equal to some composition of many morphisms) are searched and tend to be excited. 
 The searching criteria and the excitation probability is provided by (the localization of) $\mu$. 
 \item As a special case of rule 1 or 2, inversely-directed morphisms of excited morphisms tend to be excited.

 \end{itemize}
 
\item Relaxation processes (relatively slow) occur under the rule below:
\begin{itemize}
\item (Anti-Fork rule). For the pairs of morphisms sharing the domains but morphisms between their codomains 
(which may be equal to some composition of many morphisms) are not excited, is relaxed. 
\end{itemize}

\item  As a result, an indeterminate category is defined and a family of excited morphisms becomes as 
indeterminate natural transformation (=a family of subsystems of natural transformations) from (a subsystem of) 
the base-of-metaphor functor $f\backslash \mathcal{C}$ to some indeterminate functor (=metaphor functor): The INT (indeterminate natural transformation) 
provides the meaning of the metaphor.

\item When the INT above become stable (long-term surviving), comprehensive (the domain subcategory of functor defined by INT is large) and 
``valuable''(at this point out of our scheme), 
the weight of components of the INT is increased. In other words, it will give the change of weight $\mu$: The change of our world\MF{-}view.

\end{itemize}

\section{4. Application of TINT}
In this section we exemplify an application of TINT to explain a metaphor comprehension. 
We took an example of metaphor from a Japanese poem ``Tsuchi (Soil in English)''\cite{Miyoshi1932} written by Tatsuji Miyoshi.
Free translation of this poem is that:

\begin{quote}
	An ant\\
	 pull a wing of a butterfly.\\
	 O\\
	It is like yachting.
\end{quote}

For simplicity, we abstract the metaphor ``Wing is like a Sail'' from this poem to see the dynamical creation of meanings.
Figure \ref{Fig_Miyoshi_1} and \ref{Fig_Miyoshi_2} represent this metaphor according to  TINT.
Before the poem comes to our mind, 
each 
``Sail'' and ``Wing'' have made the coslice category ``$Sail \backslash \mathcal{C}$'' and ``$Wing \backslash \mathcal{C}$''.
When we read the poem ``Tsuchi'', a morphism $f$ from ``Wing'' to ``Sail'' is excited.
Then the base-of-metaphor functor 
$f\backslash \mathcal{C}:=(\cdot)\circ f: Sail\backslash \mathcal{C}\longrightarrow 
Wing\backslash \mathcal{C}$ (i.e. ``$(\cdot)$for $Sail$'' is mapped to ````$(\cdot)$for $Sail$''for $Wing$ '')
 is canonically constructed through excited $f$.
With the dynamics of excitement/relaxation, the base-of-metaphor functor 
causes the construction of the meaning of the metaphor as indeterminate natural transformation as follows:

\begin{figure}[t]
\begin{center}
\includegraphics[width = 7.5cm]{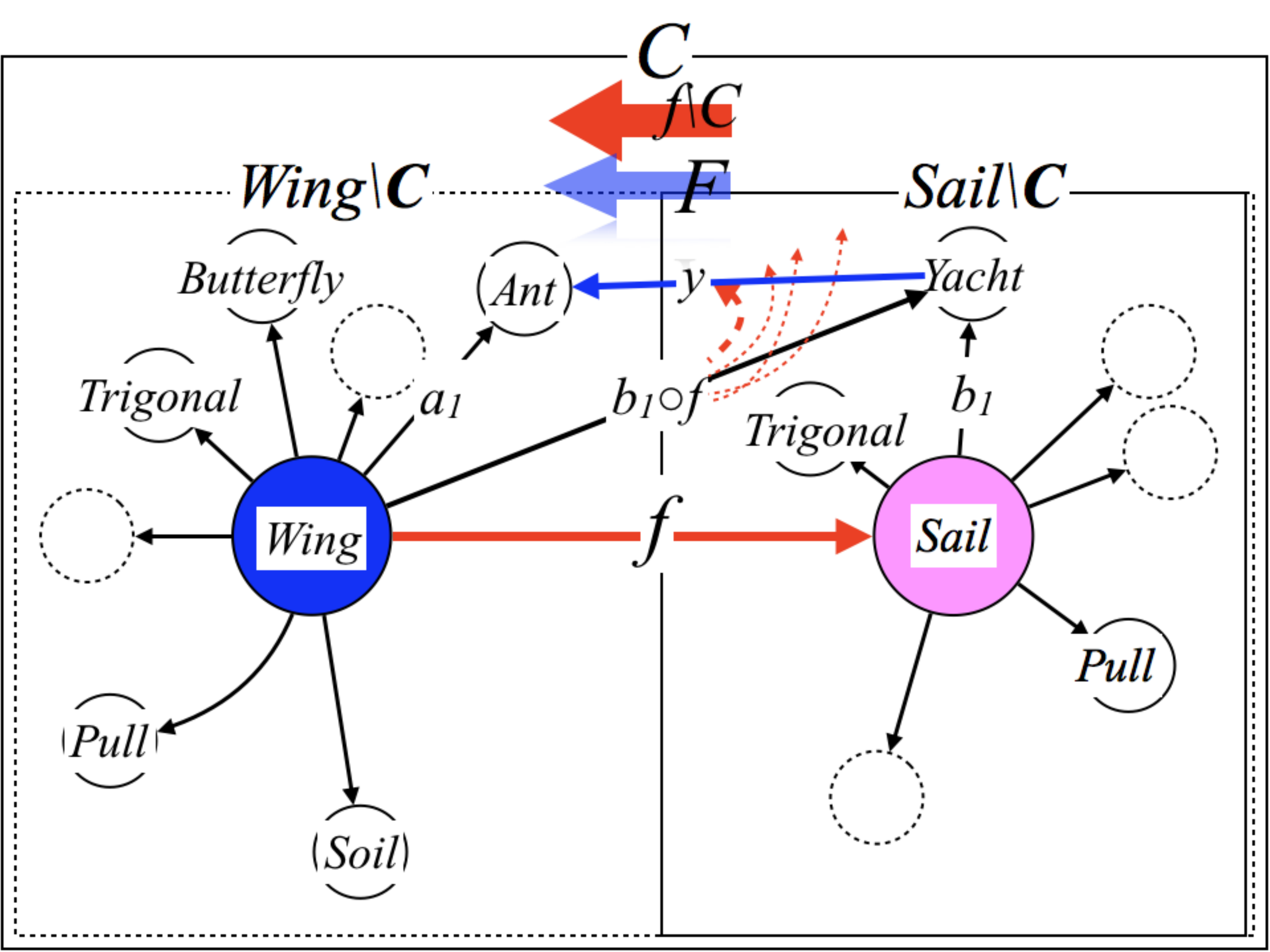}
\end{center}
\caption{\MF{``Tsuchi (Soil)'' : The base-of-metaphor functor $f \backslash \mathcal{C}$ is canonically constructed through excited f. Then the morphism from \textit{Yacht} will be searched according to the fork rule, and \textit{``Yacht to Ant''} will be excited.} }
\label{Fig_Miyoshi_1} 
\end{figure}

\begin{figure}[t]
\begin{center}
\includegraphics[width = 7.5cm]{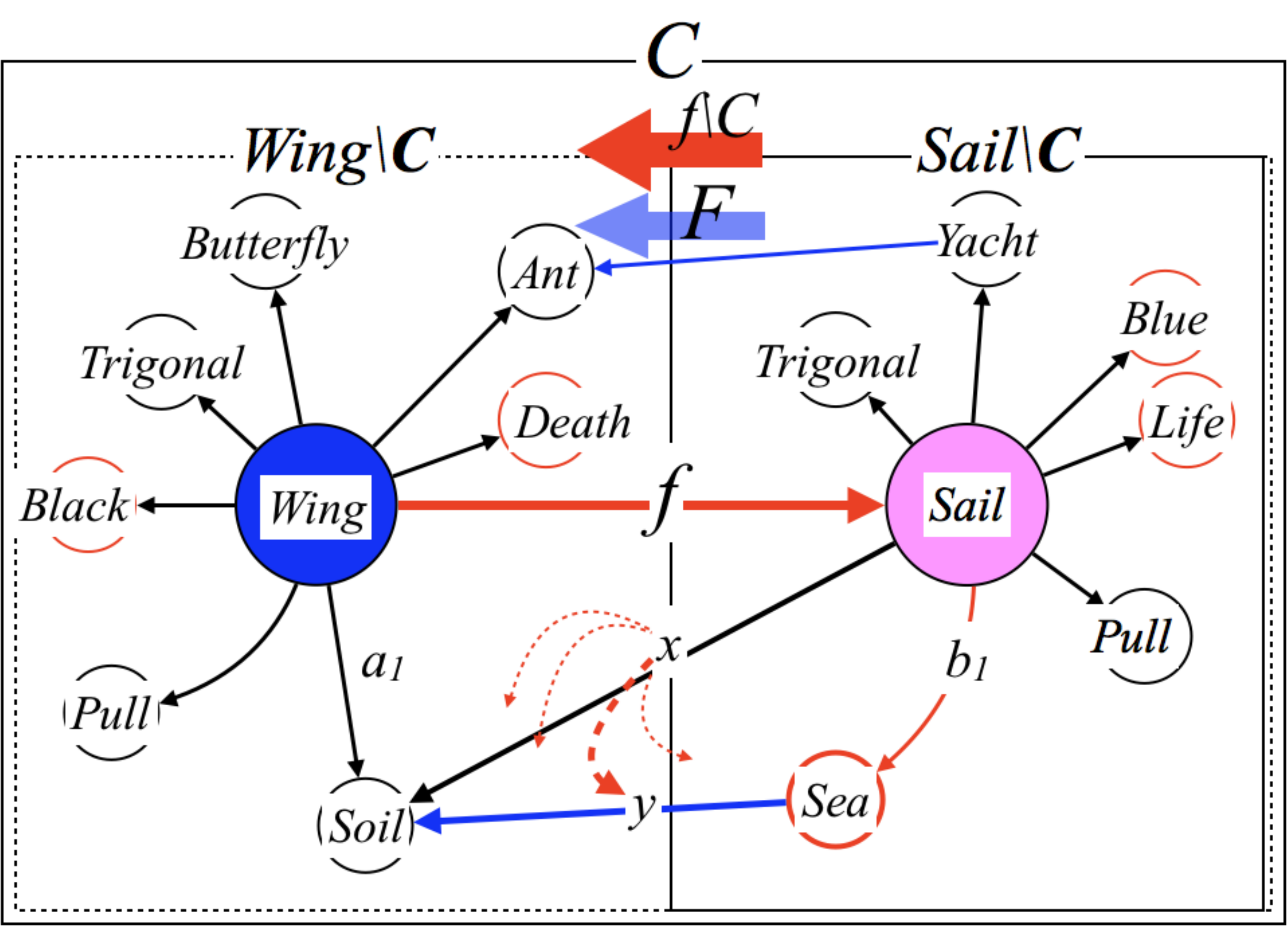}
\end{center}
\caption{\MF{``Tsuchi (Soil)'' : The ``hidden'' image of \textit{Sea} is excited according to the fork rule and probabilistic weight of $\mu$. These processes can contribute to associate other images such as ``black'', ``blue'', ``death'', ``life'', and so on.} }
\label{Fig_Miyoshi_2}
\end{figure}

For example, in figure \ref{Fig_Miyoshi_1}, according to the basic rule, ``$Wing$ to $Yacht$'' is excited.\footnote{Here we naturally assume the the morphism ``$Sail$ to $Yacht$ is already excited 
because of the expression of the poem.} By the fork rule, the morphism from $Yacht$ 
will be searched to be excited. In this context\MF{,} it is likely that the morphism ``$Yacht$ to $Ant$'' is excited. Although the weight of the morphism is not 
so much large in total $\mu$, it is natural to think that its weight will be evaluated as some local maxima since it has relatively natural 
``transit point'' like the image of ``moving on some plane entity''. The fork rule will excite many morphisms from the $Sail$ side to $Wing$ side\footnote{Since the fork rule is symmetric, the inversely directed morphisms will be also excited. Here we focus on the ``$Sail$ to $Wing$'' 
direction just for simplicity.}, 
including rather \MF{a} trivial morphism from ``$triagonal$ to $triagonal$'' etc. In general, excitement/relaxation 
dynamics make these newly constructed morphisms as the part of natural transformations $t$ from the base-of-metaphor functor to some (subsystem of) other 
functor $F$
:
\begin{align*}
\xymatrix{
\ar@<-11pt>@{}[d]_~="t1" \ar@<11pt>@{}[r]^~="t3"
& b_2 & b_1 \ar[l]_i \ar@<11pt>@{}[d]^~="t2" \\
f \backslash \mathcal{C} \ar@{=>}[d]_t & b_2 \circ f \ar[d]_{t_{b2}} & b_1 \circ f \ar[d]^{t_{b1}} \ar[l]_{i} \\
F \ar@<-11pt>@{}[r]_~="t4" & a_2 & a_1 \ar[l]^{h}
\ar@{.} "t1";"t2"
\ar@{.} "t3";"t4"
}
\end{align*}
Here $b_1, b_2$ and $i$ are in $Sail\backslash \mathcal{C}$ and $a_1,a_2$, $h$, $i$, $t_{b_1}, t_{b_2}$ are 
$Wing\backslash \mathcal{C}$.\footnote{$i$ (together with $b_1\circ f$ and $b_2 \circ f$) 
can be considered as a morphism in $Wing\backslash \mathcal{C}$. } Note that by the anti-fork rule, the family of morphism satisfying 
non-trivial commutative diagrams as above tend to survive than just those with only trivial commutative diagrams. 
As a result, these morphisms construct a family of (subsystem of) natural transformations 
from the base-of-metaphor functor $f \backslash \mathcal{C}$ to a new functor $F$ (metaphor as \MF{a} functor) : 
\MF{The} indeterminate natural transformation from the base-of metaphor functor as the creation of the meanings of metaphor.    

The meaning creation explained above also includes fairly non-trivial process. Let us focus on the pair of ``$Wing$ to $Soil$'' and $f$. 
By the fork rule, some morphism ``$Sail$ to $Soil$'' will be searched to be excited. Although it is not easy to find direct path, 
it will be searched some indirect path $x=y\circ b_1$, with the ``transit point'' $Sea$.
As a result, the ``hidden'' image of $Sea$ is excited and the morphism $y$ become a part of indeterminate natural transformation. 
In other words, the indeterminate natural transformation makes \MF{a} non-trivial fusion of the world of soil and the world of \MF{the} sea. 
Maybe these process continues to make non-trivial fusions of ``black'' and ``blue'' or ``death'' and ``life'', for example (See Figure \ref{Fig_Miyoshi_2} ). 
This continuation phenomena will expand and fertilize our meaning of this metaphor. 
In short, TINT can also explain the function of a metaphor as the method of creating ``profound'' meanings.

\section{5. Contribution of TINT}
\subsection{Characteristics of TINT}
Now we will make clear the characteristics of TINT and discuss the contributions for metaphor studies.
Characteristics of TINT can be summarized as below:

\begin{itemize}
	\item TINT represents meanings of images and metaphors as a structure of images using coslice categories.
	\item TINT represents the structure-to-structure interaction between a target and a source as 
        the construction of indeterminate natural transformation.
	\item TINT explains the process of comprehension of novel metaphor with co-existence of various interpretations by its indeterminacy.
	\item TINT models both the static meaning space by weight of morphisms and the dynamical processes of metaphor
         comprehension by excitement/relaxation dynamics.
\end{itemize}

From our viewpoint, the most important process of the metaphor 
comprehension is that the constructing novel meanings by interacting a meaning of a target and a source.
The meanings of the target and the source are determined by \MF{the} relationship among images, which can be regarded as structures of images.
Therefore, it is needed and natural to explain the interaction between the target and the source by structure-to-structure interaction.

Previous studies seem not to be able to treat this structural interaction effectively.
Although \citeA{Bowdle2005} tried to introduce structure-to-structure interaction to study metaphor,
they, eventually, did one-to-one image mapping and did not deal with structural interaction itself.
Other studies modeled the process of  metaphor comprehension as a composition of vectors of a target and a source which represented 
in a semantic space with metrics.
For example, \citeA{Kintsch2001} proposed the \textit{Predication algorithm} to model metaphor comprehension in a structural manner.
Predication algorithm represents the meaning of words as a vector in a high-dimensional semantic space using Latent Semantic Analysis.
Then, 
in short,
some vectors which are most relevant to target and source are picked up 
and combine with the centroid of the target and the source.
Although It seems very important that
predication algorithm tries to model the context of the meanings of the target and the source,
this algorithm compresses the information of contexts and interaction between meanings as  the composition of vectors despite their philosophy.
For us the vector composition seems to be too simple to represent the structural interaction, 
at least to explain the novel construction of profound meanings which we discuss in the previous section.

TINT provides the framework to deal with this interaction between the target and the source in a structural manner by virtue of the above characteristics.
To be more precise, TINT represents all essential pieces such as the meaning of image, metaphor, the process of comprehension as relationships among images
using categorical theory and stochastic processes.
That is why we believe that TINT naturally models the process of metaphor comprehension as an interaction of the structure of these relationships.

\subsection{Prospects}
The present paper proposed the fundamental concepts and working hypothesis of TINT as the first step.
In the next step, we have to sophisticate the ``indeterminate'' process as a mathematical and cognitive studies manner.
Furthermore, we should try some simulation based on TINT using some corpus, which reflects the space of images and semantic structures to check TINT's adequacy.
This line of thinking is compatible with the idea of word embedding, so we can approximately  use these semantic space with metrics which can be 
calculated by appropriate software such as word2vec\cite{mikolov2013a,mikolov2013b}.
Since the algorithm of TINT is quite simple, it will probably be implemented and simulated easily if we get appropriate semantic 
space from these corpora and software.
We also planed conducting some experiments employing subjects to verify TINT.
For example, we can control excitation/inhibition of images using priming technique and check the changes of interpretations of metaphors. 

To sum up, we prospect that TINT will work as a new platform of the interplay between the fields related to meanings, metaphors\MF{,} and morphisms.

\end{document}